\documentclass[12pt,a4paper]{amsart}
\usepackage{t1enc}
\usepackage[latin1]{inputenc}
\usepackage[english]{babel}
\usepackage{hyperref}
\usepackage{graphicx}

\usepackage{latexsym}
\usepackage{amsmath,amsfonts,amssymb,amsthm}
\usepackage{mathptmx}

\begin{document}

\title {
	Complex Gradient Systems
}

\author{Giuseppe Tomassini}
\address{
	G. Tomassini:
	Scuola Normale Superiore,
	Piazza dei Cavalieri,
	7 --- I-56126 Pisa, Italy}
\email{g.tomassini@sns.it}

\author{Sergio Venturini}
\address{
	S. Venturini:
	Dipartimento Di Matematica,
	Universit\`{a} di Bologna,
	\,\,Piazza di Porta S. Donato 5 ---I-40127 Bologna,
	Italy}
\email{venturin@dm.unibo.it}

\keywords{
	Complex manifolds, CR Manifolds, Symplectic geometry, Contact geometry}
\subjclass[2000]{Primary 32W20, 32V40 Secondary 53D10, 35Fxx}

%
%
\def\R{{\rm I\kern-.185em R}}
\def\RR{\mathbb{R}}
\def\C{{\rm\kern.37em\vrule height1.4ex width.05em depth-.011em\kern-.37em C}}
\def\CC{\mathbb{C}}
\def\N{{\rm I\kern-.185em N}}
\def\NN{\mathbb{N}}
\def\Z{{\bf Z}}
\def\ZZ{\mathbb{N}}
\def\Q{{\mathchoice{\hbox{\rm\kern.37em\vrule height1.4ex width.05em 
depth-.011em\kern-.37em Q}}{\hbox{\rm\kern.37em\vrule height1.4ex width.05em 
depth-.011em\kern-.37em Q}}{\hbox{\sevenrm\kern.37em\vrule height1.3ex 
width.05em depth-.02em\kern-.3em Q}}{\hbox{\sevenrm\kern.37em\vrule height1.3ex
width.05em depth-.02em\kern-.3em Q}}}}
\def\P{{\rm I\kern-.185em P}}
\def\H{{\rm I\kern-.185em H}}
%
\def\Aleph{\aleph_0}
\def\ALEPH#1{\aleph_{#1}}
\def\sset{\subset}\def\ssset{\sset\sset}
%
\def\bar#1{\overline{#1}}
\def\dim{\mathop{\rm dim}\nolimits}
\def\half{\textstyle{1\over2}}
\def\Half{\displaystyle{1\over2}}
\def\mlog{\mathop{\half\log}\nolimits}
\def\Mlog{\mathop{\Half\log}\nolimits}
\def\Det{\mathop{\rm Det}\nolimits}
\def\Hol{\mathop{\rm Hol}\nolimits}
\def\Aut{\mathop{\rm Aut}\nolimits}
\def\Re{\mathop{\rm Re}\nolimits}
\def\Im{\mathop{\rm Im}\nolimits}
\def\Ker{\mathop{\rm Ker}\nolimits}
\def\Fix{\mathop{\rm Fix}\nolimits}
\def\Exp{\mathop{\rm Exp}\nolimits}
\def\sp{\mathop{\rm sp}\nolimits}
\def\id{\mathop{\rm id}\nolimits}
\def\Rank{\mathop{\rm rk}\nolimits}
\def\Trace{\mathop{\rm Tr}\nolimits}
\def\Res{\mathop{\rm Res}\limits}
\def\cancel#1#2{\ooalign{$\hfil#1/\hfil$\crcr$#1#2$}}
\def\prevoid{\mathrel{\scriptstyle\bigcirc}}
\def\void{\mathord{\mathpalette\cancel{\mathrel{\scriptstyle\bigcirc}}}}
\def\n{{}|{}\!{}|{}\!{}|{}}
\def\abs#1{\left|#1\right|}
\def\norm#1{\left|\!\left|#1\right|\!\right|}
\def\nnorm#1{\left|\!\left|\!\left|#1\right|\!\right|\!\right|}
%
\def\upperint{\int^{{\displaystyle{}^*}}}
\def\lowerint{\int_{{\displaystyle{}_*}}}
\def\Upperint#1#2{\int_{#1}^{{\displaystyle{}^*}#2}}
\def\Lowerint#1#2{\int_{{\displaystyle{}_*}#1}^{#2}}
%
\def\rem #1::#2\par{\medbreak\noindent{\bf #1}\ #2\medbreak}
\def\proclaim #1::#2\par{\removelastskip\medskip\goodbreak{\bf#1:}
\ {\sl#2}\medskip\goodbreak}
\def\ass#1{{\rm(\rmnum#1)}}
\def\assertion #1:{\Acapo\llap{$(\rmnum#1)$}$\,$}
\def\Assertion #1:{\Acapo\llap{(#1)$\,$}}
\def\acapo{\hfill\break\noindent}
\def\Acapo{\hfill\break\indent}
\def\proof{\removelastskip\par\medskip\goodbreak\noindent{\it Proof.\/\ }}
\def\prova{\removelastskip\par\medskip\goodbreak
\noindent{\it Dimostrazione.\/\ }}
\def\risoluzione{\removelastskip\par\medskip\goodbreak
\noindent{\it Risoluzione.\/\ }}
\def\qed{{$\Box$}\par\smallskip}
\def\BeginItalic#1{\removelastskip\par\medskip\goodbreak
\noindent{\it #1.\/\ }}
\def\iff{if, and only if,\ }
\def\sse{se, e solo se,\ }
\def\rmnum#1{\romannumeral#1{}}
\def\Rmnum#1{\uppercase\expandafter{\romannumeral#1}{}}
\def\smallfrac#1/#2{\leavevmode\kern.1em
\raise.5ex\hbox{\the\scriptfont0 #1}\kern-.1em
/\kern-.15em\lower.25ex\hbox{\the\scriptfont0 #2}}
%
\def\Left#1{\left#1\left.}
\def\Right#1{\right.^{\llap{\sevenrm
\phantom{*}}}_{\llap{\sevenrm\phantom{*}}}\right#1}
%
%
%
\def\dimens{3em}
\def\symb[#1]{\noindent\rlap{[#1]}\hbox to \dimens{}\hangindent=\dimens}
\def\references{\bigskip\noindent{\bf References.}\bigskip}
\def\art #1 : #2 ; #3 ; #4 ; #5 ; #6. \par{#1, 
{\sl#2}, #3, {\bf#4}, (#5), #6.\par\smallskip}
\def\book #1 : #2 ; #3 ; #4. \par{#1, {\bf#2}, #3, #4.\par\smallskip}
\def\freeart #1 : #2 ; #3. \par{#1, {\sl#2}, #3.\par\smallskip}
%
%
%
%
\def\name{\hbox{Sergio Venturini}}
\def\snsaddress{\indent
\vbox{\bigskip\bigskip\bigskip
\name
\hbox{Scuola Normale Superiore}
\hbox{Piazza dei Cavalieri, 7}
\hbox{56126 Pisa (ITALY)}
\hbox{FAX 050/563513}}}
\def\cassinoaddress{\indent
\vbox{\bigskip\bigskip\bigskip
\name
\hbox{Universit\`a di Cassino}
\hbox{via Zamosch 43}
\hbox{03043 Cassino (FR)}
\hbox{ITALY}}}
\def\bolognaaddress{\indent
\vbox{\bigskip\bigskip\bigskip
\name
\hbox{Dipartimento di Matematica}
\hbox{Universit\`a di Bologna}
\hbox{Piazza di Porta S. Donato 5}
\hbox{40127 Bologna (BO)}
\hbox{ITALY}
\hbox{venturin@dm.unibo.it}
}}
\def\homeaddress{\indent
\vbox{\bigskip\bigskip\bigskip
\name
\hbox{via Garibaldi, 7}
\hbox{56124 Pisa (ITALY)}}}
\def\doubleaddress{
\vbox{
\hbox{\name}
\hbox{Universit\`a di Cassino}
\hbox{via Zamosch 43}
\hbox{03043 Cassino (FR)}
\hbox{ITALY}
\smallskip
\hbox{and}
\smallskip
\hbox{Scuola Normale Superiore}
\hbox{Piazza dei Cavalieri, 7}
\hbox{56126 Pisa (ITALY)}
\hbox{FAX 050/563513}}}
\def\sergio{{\rm\bigskip
\centerline{Sergio Venturini}
\leftline{\bolognaaddress}
\bigskip}}
%
%
%
%
%
%

\newtheorem{theorem}{Theorem}[section]
\newtheorem{proposition}{Proposition}[section]
\newtheorem{lemma}{Lemma}[section]
\newtheorem{corollary}{Corollary}[section]
\newtheorem{remark}{Remark}[section]
\newtheorem{definition}{Definition}[section]

\newtheorem{teorema}{Teorema}[section]
\newtheorem{proposizione}{Proposizione}[section]
\newtheorem{corollario}{Corollario}[section]
\newtheorem{osservazione}{Osservazione}[section]
\newtheorem{definizione}{Definizione}[section]
\newtheorem{esempio}{Esempio}[section]
\newtheorem{esercizio}{Esercizio}[section]
\newtheorem{congettura}{Congettura}[section]

\bibliographystyle{abbrv}


\def\CRMan{M}
\def\CXMan{\widetilde{M}}
\def\CXdim{{n+k}}
\def\CRdim{n}
\def\CRcodim{k}
\def\Tan{T}
\def\CTan{T_\mathbb{C}}
\def\HTan{T^{(1,0)}_\mathbb{C}}
\def\pMan{p}
\def\CXJ{J}

\def\CR{{\rm CR}}
\def\dd{{\rm d}}
\def\dc{{\dd^c}}

\def\CGSBase{\mathcal{V}}
\def\CGSBaseC{\CGSBase^{\mathbb{C}}}
\def\CGSRep{\rho}
\def\CGSCRep{\CGSRep^{\mathbb{C}}}
\def\CGSGradmap{U}
\def\CGSgradmap{u}
\def\CGSCRTan{{\mathcal{H}}}
\def\CGSVect{V}
\def\CGSVectB{W}
\def\CGSDistr{{\mathcal{D}}}
\def\CGSRDistr#1{{\mathcal{D}^{\mathbb{R}}_{#1}}}
\def\CGSCDistr#1{{\mathcal{D}^{\mathbb{C}}_{#1}}}
\def\CGSExp{G}
\def\CGSCExp{{\widetilde\CGSExp}}
\def\CGSIExp{F}
\def\FIsoDom{{\widetilde\CGSCexpDom}}
\def\CGSCexpDom{D}
\def\CGSCDom{N}
\def\ucomp{u}
\def\CGSVBase{\xi}
\def\CGSGroup{G}
\def\CGSCGroup{G^{\mathbb{C}}}
\def\CGSAlg{\mathfrak{g}}
\def\CGSAlgD{\mathfrak{g}^*}
\def\CGSCAlg{\mathfrak{g}^{\mathbb{C}}}
\def\CGSGElem{g}
\def\CGSGElemB{g_1}
\def\CGSAElem{V}
\def\CGSAElemB{W}
\def\CGSGGMap#1{\CGSGradmap_{#1}}
\def\CGSGRep#1{\CGSRep_{#1}}

\def\ixa{\alpha}
\def\ixb{\beta}
\def\ixc{\gamma}
\def\ixd{\delta}
\def\iya{\mu}

\def\VField{X}

\def\res{\mathop{\hbox{\vrule height 7pt width .5pt depth 0pt\vrule height .5pt width 6pt depth 0pt}}\nolimits}

\def\LieBracket#1#2{\left[{#1},{#2}\right]}
\def\der[#1]{\frac{\partial}{\partial {#1}}}
\def\pder#1#2{\frac{\partial {#1}}{\partial {#2}}}

\begin{abstract}
Let $\CXMan$ be a complex manifold of complex dimension $\CXdim$.
We say that the functions $\ucomp_1,\ldots,\ucomp_\CRcodim$
and the vector fields $\CGSVBase_1,\ldots,\CGSVBase_\CRcodim$
on $\CXMan$ form a \emph{complex gradient system} if
$\CGSVBase_1,\ldots,\CGSVBase_\CRcodim,\CXJ\CGSVBase_1,\ldots,\CXJ\CGSVBase_\CRcodim$
are linearly independent at each point $\pMan\in\CXMan$ and generate
an integrable distribution of $\Tan\CXMan$ of dimension $2\CRcodim$
and $\dd\ucomp_\ixa(\CGSVBase_\ixb)=0$, $\dc\ucomp_\ixa(\CGSVBase_\ixb)=\delta_{\ixa\ixb}$
for $\ixa,\ixb=1,\ldots,\CRcodim$.

We prove a Cauchy theorem for such complex gradient systems with
initial data along a $\CR-$submanifold of type $(\CRdim,\CRcodim)$.

We also give a complete local characterization for
the complex gradient systems which are
\emph{holomorphic} and \emph{abelian},
which means that the vector fields
$\CGSVBase_\ixa^c=\CGSVBase_\ixa-i\CXJ\CGSVBase_\ixa$,
$\ixa=1,\ldots,\CRcodim$ are holomorphic and satisfy
$\big[\CGSVBase_\ixa^c,\bar{\CGSVBase_\ixb^c}\bigr]=0$
for each $\ixa,\ixb=1,\ldots,\CRcodim$.

\end{abstract}
\maketitle



\section{\label{section:Introduction}Introduction}
{
Let $\CXMan$ be a complex manifold, $\Tan\CXMan$ its (real) tangent
space endowed by its complex structure $\CXJ$.

In \cite{article:TV2ArxivContactMA} the authors
introduced a geometric tool named
\emph{one dimensional calibrated foliation}
on the complex manifold $\CXMan$. It consists of a real function $\ucomp:\CXMan\to\R$
and a vector field $\CGSVBase\in\Gamma(\CXMan,\Tan\CXMan)$
which satisfy the conditions
\begin{eqnarray*}
	[\CGSVBase,\CXJ\CGSVBase]=0,\\
	\dd\ucomp(\CGSVBase)=0,\\
	\dc\ucomp(\CGSVBase)=1.
\end{eqnarray*}
Here $[\cdot,\cdot]$ is the Poisson Lie bracket between vector fields
and $\dc\ucomp(\CGSVBase)=-\dd\ucomp\bigl(\CXJ(\CGSVBase)\bigr)$.

Among the results proved in \cite{article:TV2ArxivContactMA} 
there is a Cauchy-like theorem for one dimensional calibrated foliation (see Theorem 3.1) which states the following:
if $\CRMan\sset\CXMan$ is a real hypersurface and
$\CGSVBase_0$ is a vector field on $\CRMan$ which is transversal to the 
holomorphic tangent space to $\CRMan$,
then, under the assumption that the integral curves $t\mapsto\gamma(t)$
of $\CGSVBase_0$ are real analytic, there exists a (unique) one dimensional calibrated foliation
$(\CGSVBase,\ucomp)$, defined in a suitable neighbourhood of $\CRMan$
in $\CXMan$, such that the vector field $\CGSVBase$ extends $\CGSVBase_0$.

The notion one dimensional calibrated foliation was motivated by the problem of finding for $\CRMan$ an equation satisfying the homogeneous complex Monge-Amp\`ere equation. The method has subsequently been applied in \cite{article:TV3AdaptedTube},
\cite{article:TV3ArxivAdaptedTube} to prove the existence of adapted complex structures on the symplectization of a pseudo-Hermitian manifolds.
The key point in proving the existence of a calibrated foliation
is the construction of a function $\CGSCExp(z,\pMan):\CGSCexpDom\to\CXMan$,
where $\CGSCexpDom\sset\CRMan\times\C$ is an open neighbourhood
of $\CRMan\times\{0\}$, which is holomorphic in $z$ for each $\pMan\in\CRMan$ and for $z=t$ real the map $t\mapsto\CGSCExp(t,\pMan)$
is the integral curve of the vector field $\CGSVBase_0$ such that
$\CGSCExp(0,\pMan)=\pMan$. This idea goes back to Duchamp and Kalka (see \cite{article:DuchampKalka}).  

The purpose of this paper is to provide a natural higher dimensional generalization of the notion of one dimensional calibrated foliations.

Let $\CXMan$ be a complex manifold of complex dimension $\CXdim$.
We say that the functions $\ucomp_1,\ldots,\ucomp_\CRcodim$
and the vector fields $\CGSVBase_1,\ldots,\CGSVBase_\CRcodim$
on $\CXMan$ form a \emph{complex gradient system}
(of dimension $\CRcodim$) if
$$
\CGSVBase_1,\ldots,\CGSVBase_\CRcodim,\CXJ(\CGSVBase_1),\ldots,\CXJ(\CGSVBase_\CRcodim)
$$
are linearly independent at each point $\pMan\in\CXMan$ and generate
an integrable distribution of $\Tan\CXMan$ of dimension $2\CRcodim$
and 
$$
\dd\ucomp_\ixa(\CGSVBase_\ixb)=0, \>\>\dc\ucomp_\ixa(\CGSVBase_\ixb)=\delta_{\ixa\ixb}
$$
for $\ixa,\ixb=1,\ldots,\CRcodim$.
Here $\delta_{\ixa\ixb}$ is the usual Kronecker symbol.

In a more intrinsic way a complex gradient system is given by
a real vector space $\CGSBase$ of dimension $\CRcodim$,
a linear monomorphism $\CGSRep:\CGSBase\to\Gamma(\CXMan,\Tan\CXMan)$,
the \emph{representation} map, and a map $\CGSGradmap:\CXMan\to\CGSBase$,
the \emph{gradient map}, which satisfy
\begin{eqnarray*}
	\dd\CGSGradmap\bigl(\CGSRep(\CGSVect)\bigr)=0,\\
	\dc\CGSGradmap\bigl(\CGSRep(\CGSVect)\bigr)=\CGSVect
\end{eqnarray*}
for each $\CGSVect\in\CGSBase$.

The name ``complex gradient system''
(instead of ``calibrated foliation'')
arises from the fact that there are examples of triples
$(\CGSBase,\CGSRep,\CGSGradmap)$  where
$\CGSBase=\CGSAlg$ is the Lie algebra of a Lie group $\CGSGroup$
which is a compact real form of a reductive complex Lie group $\CGSCGroup$
and the map $\CGSGradmap$
(with the identification of $\CGSAlg=\CGSAlgD$ with its dual $\CGSAlgD$ by the Killing form)
is the moment map associated to a symplectic action of $\CGSGroup$.
It is customary in the symplectic geometry to call such a moment map
as ``gradient map'' (see e.g. \cite{article:HeinznerSchuetzdellerConvexGradMap})
and hence the name ``complex gradient system''.

For a general reference on symplectic geometry and moment maps theory see e.g.
\cite{book:DusaSalamonIntro} and \cite{book:AnaCannasDaSilva}.
Basic definitions and notions in $\CR-$geometry can be found in \cite{book:DragominTomassini}.

Let now describe in more details the content of the paper.

In Section \ref{section:CGS} contains the elementary
properties of a complex gradient system. In particular,
any complex gradient system satisfies the formal
commutativity property
\begin{eqnarray*}
	\bigl[\CGSCRep(\CGSVect),\CGSCRep(\CGSVectB)]&=&0
\end{eqnarray*}
for each $\CGSVect,\CGSVectB\in\CGSBase$,
where 
$\CGSCRep(\CGSVect)=\frac{1}{2}
	\left(\CGSRep(\CGSVect)-i\CXJ\CGSRep(\CGSVect)\right)
$
(and similarly
$\CGSCRep(,\CGSVectB)=\frac{1}{2}
	\left(\CGSRep(\CGSVectB)-i\CXJ\CGSRep(\CGSVectB)\right)
$)
is the complex vector field of type $(1,0)$ naturally associated
to the real vector field $\CGSRep(\CGSVect)$ (resp. $\CGSRep(\CGSVectB)$).
See Theorem \ref{stm::CGSCommuting}.

In section \ref{section:Cauchy} we solve a Cauchy problem for a
complex gradient system on a complex manifold $\CXMan$ of dimension
$\CXdim$ with initial data on a $\CR-$submanifold 
of $\CXMan$ of type $(\CRdim,\CRcodim)$ (Theorem \ref{stm::CGSCauchy}).

In section \ref{section:ExampleLieGroups} we give a couple
of examples applying our construction to the case of the
complexification of a real Lie group $\CGSGroup$. In particular,
we find explicitly the complex gradient system associated
to the standard representation of the Lie algebra
$\CGSAlg$ of $\CGSGroup$ as left invariant vector fields
on $\CGSGroup$.

Finally, in the last section we give a complete local description of
any \emph{abelian holomorphic} complex gradient system
$(\CGSBase,\CGSRep,\CGSGradmap)$, where
abelian means that
$\bigl[\CGSCRep(\CGSVect),\bar{\CGSCRep(\CGSVectB)}\bigr]=0$
for each pair of vectors $\CGSVect,\CGSVectB\in\CGSBase$
and holomorphic means that $\CGSCRep(\CGSVect)$ is
a holomorphic vector field on $\CXMan$ for each $\CGSVect\in\CGSBase$.
See Theorem \ref{stm::HoloAbelian} for details.

}

\section{\label{section:CGS}Complex gradient  systems}
{ 
Let $\CXMan$ be a complex manifold of complex dimension $\CXdim$,
$\Tan\CXMan$ its (real) tangent 	space endowed by its complex structure $\CXJ$.
Let $\CGSBase$ be a real vector space and let
$\CGSRep:\CGSBase\to\Gamma(\CXMan,T\CXMan)$ be a linear map.
We denote by  $\CGSRDistr\CGSRep\sset\Tan\CXMan$ the distribution
generated by the vector fields of the form $\CGSRep(\CGSVect)$, $\CGSVect\in\CGSBase$.
We also denote by $\CGSCDistr\CGSRep\sset\Tan\CXMan$ the distribution
\begin{eqnarray*}
	&&\CGSCDistr\CGSRep=\CGSRDistr\CGSRep+J\bigl(\CGSRDistr\CGSRep\bigr)
\end{eqnarray*}
generated by the vector fields of the form
$\CGSRep(\CGSVect)$ and $J\bigl(\CGSRep(\CGSVect)\bigr)$, $\CGSVect\in\CGSBase$.

\begin{definition}
A \emph{complex gradient system} of dimension $\CRcodim$ on $\CXMan$ is a triple
\begin{eqnarray*}
	(\CGSBase,\CGSRep,\CGSGradmap)
\end{eqnarray*}
where:
\begin{enumerate}
\item $\CGSBase$ is a real vector space of dimension $\CRcodim$;
\item $\CGSRep:\CGSBase\to\Gamma(\CXMan,T\CXMan)$ is a $\R-$linear map;
\item $\CGSGradmap:\CXMan\to\CGSBase$ is a smooth map.
\end{enumerate}
which satisfies
\begin{enumerate}
\item[i)] for each $\CGSVect\in\CGSBase$ the vector field $\CGSRep(\CGSVect)$
is smooth and we have the identities
\begin{eqnarray*}
	&&\dd\CGSGradmap\bigl(\CGSRep(\CGSVect)\bigr)=0,\\
	&&\dc\CGSGradmap\bigl(\CGSRep(\CGSVect)\bigr)=\CGSVect.
\end{eqnarray*}
\item[ii)] the distribution $\CGSCDistr\CGSRep\sset\Tan\CXMan$ is integrable.
\end{enumerate}
\end{definition}

The maps $\CGSRep$ and $\CGSGradmap$ are said respectively
the \emph{representation} and the \emph{gradient map}
of the complex gradient system $(\CGSBase,\CGSRep,\CGSGradmap)$.

If $\{\CGSVect_1,\ldots,\CGSVect_\CRcodim\}$ is a basis of $\CGSBase$
we set
$$\CGSVBase_1=\CGSRep(\CGSVect_1),\ldots,
\CGSVBase_\CRcodim=\CGSRep(\CGSVect_\CRcodim)$$
and for some smooth functions $\CGSgradmap_1,\ldots,\CGSgradmap_\CRcodim$
we have
$$\CGSGradmap=\CGSgradmap_1\CGSVect_1+\cdot+
\CGSgradmap_\CRcodim\CGSVect_\CRcodim.$$
Then $(\CGSBase,\CGSRep,\CGSGradmap)$ is a complex gradient system if,
and only if,
\begin{eqnarray*}
	&&\dd\CGSgradmap_\ixa(\CGSVBase_\ixb)=0,
		\quad\ixa,\ixb=1,\ldots,\CRcodim,\\
	&&\dc\CGSgradmap_\ixa(\CGSVBase_\ixb)=\delta_{\ixa\ixb},
		\quad\ixa,\ixb=1,\ldots,\CRcodim,\\
\end{eqnarray*}
and $\bigl\{\CGSVBase_1,\CXJ\CGSVBase_1,\ldots,\CGSVBase_\CRcodim,\CXJ\CGSVBase_\CRcodim\bigr\}$
is a basis of an integrable distribution.

\begin{proposition}
Let $(\CGSBase,\CGSRep,\CGSGradmap)$ be a complex gradient system
on the manifold $\CXMan$.
Then
\begin{eqnarray*}
	&&\CGSRDistr\CGSRep=\CGSCDistr\CGSRep\cap\ker\dd\CGSGradmap,\\
	&&\CGSCDistr\CGSRep=\CGSRDistr\CGSRep\oplus\CXJ\CGSRDistr\CGSRep,\\
	&&\Tan\CXMan=\CXJ\CGSRDistr\CGSRep\oplus\ker\dd\CGSGradmap.
\end{eqnarray*}
\end{proposition}

{\proof
\def\vtan{v}
\def\vtanb{w}
Let $\pMan\in\CXMan$ and $\vtan\in\Tan_\pMan\CXMan$.
Assume that $\vtan\in\CGSCDistr\CGSRep$.
Then there are $\CGSVect,\CGSVectB\in\CGSBase$ such that
$\vtan=\CGSRep(\CGSVect)_\pMan+\CXJ\CGSRep(\CGSVectB)_\pMan$.
It follows
$$\dd\CGSGradmap(\vtan)=
\dd\CGSGradmap\bigl(\CGSRep(\CGSVect)_\pMan\bigr)
+\dd\CGSGradmap\bigl(\CXJ\CGSRep(\CGSVectB)_\pMan\bigr)
=-\CGSVectB,$$
whence
$$
	\vtan\in\ker\dd\CGSGradmap
	\Longleftrightarrow\CGSVectB=0
	\Longleftrightarrow\vtan=\CGSRep(\CGSVect)_\pMan\in\CGSRDistr\CGSRep.
$$
This proves the first assertion of the proposition.

As for the second one it suffices to prove that
$\CGSRDistr\CGSRep\cap\CXJ\CGSRDistr\CGSRep=0$.
Let $\vtan\in\CGSRDistr\CGSRep\cap\CXJ\CGSRDistr\CGSRep$.
Then $\vtan=\CGSRep(\CGSVect)_\pMan=\CXJ\CGSVectB_\pMan$
for some $\CGSVect,\CGSVectB\in\CGSBase$.
We then have
$$
	0=\dd\CGSGradmap\bigl(\CGSRep(\CGSVect)_\pMan\bigr)
	=\dd\CGSGradmap\bigl(\CXJ\CGSRep(\CGSVectB)_\pMan)
	=-\CGSVectB
$$
and hence $\vtan=\CXJ\CGSRep(\CGSVectB)_\pMan=0$,
as required.

Let now $\vtan\in\Tan_\pMan\CXMan$ be arbitrary and
set $\CGSVect=\dd\CGSGradmap(\vtan)$, $\vtanb=\CGSRep(\CGSVect)_\pMan$.
Clearly, $\vtan=\bigl(\vtan-\CXJ\vtanb\bigr)+\CXJ\vtanb$.
Observe that
$$
\dd\CGSGradmap\bigl(\vtan-\CXJ\vtanb\bigr)=
\dd\CGSGradmap(\vtan)-\dd\CGSGradmap\bigl(\CXJ\vtanb\bigr)=
\CGSVect-\dd\CGSGradmap\bigl(\CXJ\CGSRep(\CGSVect)_\pMan\bigr)=
\CGSVect-\CGSVect=0,
$$
i.e. $\vtan-\CXJ\vtanb)in\ker\dd\CGSGradmap$
and $\CXJ\vtanb\in\CXJ\CGSRDistr\CGSRep$.

If
$\vtan\in\CXJ\CGSRDistr\CGSRep\cap\ker\dd\CGSGradmap$,
then $\vtan=\CXJ\CGSRep(\CGSVectB)_\pMan$
for some $\CGSVectB\in\CGSBase$.
It follows that
$$
	0=\dd\CGSGradmap(\vtan)
	=\dd\CGSGradmap\bigl(\CXJ\CGSRep(\CGSVectB)_\pMan\bigr)
	=-\CGSVectB
$$
and hence $\vtan=\CXJ\CGSRep(\CGSVectB)_\pMan=\CXJ\CGSRep(0)_\pMan=0$.
This proves the last assertion of the proposition.

\qed}

\begin{definition}
Given a complex gradient system $(\CGSBase,\CGSRep,\CGSGradmap)$ we
denote $\CGSCRTan_\CGSRep$ the distribution
$\ker\,\dd\CGSGradmap\cap\ker\,\dc\CGSGradmap$.
\end{definition}

\begin{proposition}
Let $(\CGSBase,\CGSRep,\CGSGradmap)$ be a complex gradient system
on the manifold $\CXMan$.
Then
\begin{eqnarray*}
	&&\ker\,\dd\CGSGradmap=\CGSRDistr\CGSRep\oplus\CGSCRTan_\CGSRep,\\
	&&\Tan\CXMan=\CGSRDistr\CGSRep\oplus\CXJ\CGSRDistr\CGSRep\oplus\CGSCRTan_\CGSRep.
\end{eqnarray*}
\end{proposition}

{\proof
\def\vtan{v}
\def\vtanb{w}
The second equality easily follows from the first in view of the equality
$\Tan\CXMan=\CXJ\CGSRDistr\CGSRep\oplus\ker\dd\CGSGradmap$
proved in the last proposition. So it suffices to prove that
$\ker\,\dd\CGSGradmap=\CGSRDistr\CGSRep\oplus\CGSCRTan_\CGSRep$.

By definition we have $\CGSRDistr\CGSRep\sset\ker\,\dd\CGSGradmap$ and
by construction $\CGSCRTan_\CGSRep\sset\ker\,\dd\CGSGradmap$
so that $\CGSRDistr+\CGSCRTan_\CGSRep\sset\ker\,\dd\CGSGradmap$.

Let $\vtan\in\ker\,\dd\CGSGradmap$ be arbitrary.
Set $\CGSVect=\dc\CGSGradmap(\vtan)$ and
$\vtanb=\CGSRep(\CGSVect)_\pMan$.
Then we have immediately $\vtan=(\vtan-\vtanb)+\vtanb$ and
$\vtanb,\vtan-\vtanb\in\ker\,\dd\CGSGradmap$.
Observe that
$$
\dc\CGSGradmap(\vtan-\vtanb)=
\dc\CGSGradmap(\vtan)-\dc\CGSGradmap(\vtanb)=
\CGSVect-\dc\CGSGradmap\bigl(\CGSRep(\CGSVect)_\pMan\bigr)=
\CGSVect-\CGSVect=0,
$$
that is $\vtan-\vtanb\in\ker\dc\CGSGradmap$
and $\vtanb\in\CGSRDistr\CGSRep$.

Assume now that
$\vtan\in\CGSRDistr\CGSRep\cap\ker\dc\CGSGradmap$.
Then $\vtan=\CGSRep(\CGSVectB)_\pMan$
for some $\CGSVectB\in\CGSBase$.
It follows that
$$
	0=\dc\CGSGradmap(\vtan)
	=\dc\CGSGradmap\bigl(\CGSRep(\CGSVectB)_\pMan\bigr)
	=\CGSVectB
$$
and hence $\vtan=\CGSRep(\CGSVectB)_\pMan=\CGSRep(0)_\pMan=0$.
The proof is now complete.

\qed}

\begin{proposition}
Let $(\CGSBase,\CGSRep,\CGSGradmap)$ be a complex gradient system
on the manifold $\CXMan$ and $\CGSVect,\CGSVectB\in\CGSBase$.
Then
\begin{eqnarray*}
	&&\bigl[\CGSRep(\CGSVect),\CGSRep(\CGSVectB)\bigr],
	\bigl[\CGSRep(\CGSVect),J\CGSRep(\CGSVectB)\bigr],
	\bigl[J\CGSRep(\CGSVect),J\CGSRep(\CGSVectB)\bigr]
		\in\Gamma(\CXMan,\CGSRDistr\CGSRep).
\end{eqnarray*}

Moreover
\begin{eqnarray*}
	\begin{array}{lcl}
	\dd\dc\CGSGradmap\bigl(\CGSRep(\CGSVect),\CGSRep(\CGSVectB)\bigr)
			&=&-\dc\CGSGradmap\bigl([\CGSRep(\CGSVect),\CGSRep(\CGSVectB)]\bigr)\\
		\dd\dc\CGSGradmap\bigl(\CXJ\CGSRep(\CGSVect),\CXJ\CGSRep(\CGSVectB)\bigr)
			&=&-\dc\CGSGradmap\bigl([\CGSRep(\CGSVect),\CGSRep(\CGSVectB)]\bigr)\\
		\dd\dc\CGSGradmap\bigl(\CGSRep(\CGSVect),\CXJ\CGSRep(\CGSVectB)\bigr)
			&=&-\dc\CGSGradmap\bigl([\CGSRep(\CGSVect),\CXJ\CGSRep(\CGSVectB)]\bigr).
	\end{array}
\end{eqnarray*}	
\end{proposition}

{\proof
\def\VAuxA{{\xi_1}}
\def\VAuxB{{\xi_2}}
Let $\VAuxA$ be either $\CGSRep(\CGSVect)$ or $\CXJ\CGSRep(\CGSVect)$
and $\VAuxB$ be either $\CGSRep(\CGSVectB)$ or $\CXJ\CGSRep(\CGSVectB)$.
Then $\VAuxA(\CGSGradmap)$ and $\VAuxB(\CGSGradmap)$ are constant functions and hence
$$\VAuxB\bigl(\VAuxA(\CGSGradmap)\bigr)=\VAuxA\bigl(\VAuxB(\CGSGradmap)\bigr)=0.$$
This easily implies that $[\VAuxA,\VAuxB]\in\Gamma(\CXMan,\ker\dd\CGSGradmap)$.
By definition of complex gradient system,  $\CGSCDistr\CGSRep$ is an integrable distribution,
so, in view of the last proposition, we have
$\CGSRDistr\CGSRep=\CGSCDistr\CGSRep\cap\ker\dd\CGSGradmap$
and hence $[\VAuxA,\VAuxB]\in\Gamma(\CXMan,\CGSRDistr\CGSRep)$.

Using again the equalities
$\VAuxB\bigl(\VAuxA(\CGSGradmap)\bigr)=\VAuxA\bigl(\VAuxB(\CGSGradmap)\bigr)=0$
we also obtain
\begin{eqnarray*}
	\dd\dc\CGSGradmap(\VAuxA,\VAuxB)&=&
	\VAuxA\bigl(\VAuxB(\CGSGradmap)\bigr)-\VAuxB\bigl(\VAuxA(\CGSGradmap)\bigr)
	-\dd\dc\CGSGradmap\bigl([\VAuxA,\VAuxB]\bigr)\\
	&=&-\dd\dc\CGSGradmap\bigl([\VAuxA,\VAuxB]\bigr).
\end{eqnarray*}
This completes the proof of the proposition.

\qed}

\begin{corollary}
Let $(\CGSBase,\CGSRep,\CGSGradmap)$ be a complex gradient system
on the manifold $\CXMan$ and $\CGSVect,\CGSVectB\in\CGSBase$.
Then
\begin{eqnarray*}
	\begin{array}{lcl}
	\CGSRep\bigl(\dd\dc\CGSGradmap\bigl(\CGSRep(\CGSVect),\CGSRep(\CGSVectB)\bigr)\bigr)
			&=&-[\CGSRep(\CGSVect),\CGSRep(\CGSVectB)]\\
	\CGSRep\bigl(\dd\dc\CGSGradmap\bigl(\CXJ\CGSRep(\CGSVect),\CXJ\CGSRep(\CGSVectB)\bigr)\bigr)
			&=&-[\CGSRep(\CGSVect),\CGSRep(\CGSVectB)]\\
	\CGSRep\bigl(\dd\dc\CGSGradmap\bigl(\CGSRep(\CGSVect),\CXJ\CGSRep(\CGSVectB)\bigr)\bigr)
			&=&-[\CGSRep(\CGSVect),\CXJ\CGSRep(\CGSVectB)].
	\end{array}
\end{eqnarray*}	
\end{corollary}

\proof
Apply $\CGSRep$ to both sides of the last three equalities of the previous proposition
and use the identity $\dc\CGSGradmap\bigl(\CGSRep(\CGSVect)\bigr)=\CGSVect$.

\qed

\begin{corollary}
Let $(\CGSBase,\CGSRep,\CGSGradmap)$ be a complex gradient system
on the manifold $\CXMan$.
Then the distribution $\CGSRDistr\CGSRep\sset\Tan\CXMan$ is integrable.
\end{corollary}

\begin{corollary}
Let $(\CGSBase,\CGSRep,\CGSGradmap)$ be a complex gradient system
of dimension $\CRcodim$ on the complex manifold $\CXMan$ of dimension $\CXdim$.
For every $\CGSVect\in\CGSBase$ the level set $\CGSGradmap^{-1}(V)$ of the smooth function $\CGSGradmap$ is either empty or it is a $\CR-$submanifold of type $(\CRdim,\CRcodim)$.
\end{corollary}

\begin{proposition}
Let $(\CGSBase,\CGSRep,\CGSGradmap)$ be a complex gradient system
on the manifold $\CXMan$. For every $\CGSVect,\CGSVectB\in\CGSBase$ one has
\begin{eqnarray*}
	\begin{array}{lcl}
		\bigl[J\CGSRep(\CGSVect),J\CGSRep(\CGSVectB)\bigr]
			&=&\bigl[\CGSRep(\CGSVect),\CGSRep(\CGSVectB)\bigr],\\
		\bigl[J\CGSRep(\CGSVect),\CGSRep(\CGSVectB)\bigr]
			&=&-\bigl[\CGSRep(\CGSVect),J\CGSRep(\CGSVectB)\bigr].\\
	\end{array}
\end{eqnarray*}
\end{proposition}

{\proof
Since $\CXMan$ is a complex manifold the complex structure $\CXJ$ is
integrable, hence
\begin{equation}
	\CXJ\bigl[\CGSRep(\CGSVect),\CGSRep(\CGSVectB)\bigr]
		-\CXJ\bigl[J\CGSRep(\CGSVect),J\CGSRep(\CGSVectB)\bigr]=
		\bigl[J\CGSRep(\CGSVect),\CGSRep(\CGSVectB)\bigr]
		+\bigl[\CGSRep(\CGSVect),J\CGSRep(\CGSVectB)\bigr].
		\nonumber
\end{equation}
The right side of such equality belongs to
$\Gamma\bigl(\CXMan,\CXJ\CGSRDistr\CGSRep\bigr)$
while the second belongs to $\Gamma\bigl(\CXMan,\CGSRDistr\CGSRep\bigr)$.
Since $\CXJ\CGSRDistr\CGSRep\cap\CGSRDistr\CGSRep=0$ it follows that
\begin{eqnarray*}
	&&\CXJ\bigl[\CGSRep(\CGSVect),\CGSRep(\CGSVectB)\bigr]
	-\CXJ\bigl[J\CGSRep(\CGSVect),J\CGSRep(\CGSVectB)\bigr]=0,\\
	&&\bigl[J\CGSRep(\CGSVect),\CGSRep(\CGSVectB)\bigr]
	+\bigl[\CGSRep(\CGSVect),J\CGSRep(\CGSVectB)\bigr]=0
\end{eqnarray*}
and the assertion follows.

\qed}

Let $\CTan\CXMan=\C\otimes_\R\Tan\CXMan$ denote the
complexification of $\Tan\CXMan$ and $\HTan\CXMan$
the subbundle of the tangent vector of type $(1,0)$.

\begin{definition}
Let $(\CGSBase,\CGSRep,\CGSGradmap)$ be a complex gradient system.
The \emph{complexified representation}
\begin{eqnarray*}
&&\CGSCRep:\CGSBase\to\Gamma(\CXMan,\HTan\CXMan)
\end{eqnarray*}
is defined for each $\CGSVect\in\CGSBase$ by
\begin{eqnarray*}
&&\CGSCRep(\CGSVect)=\frac{1}{2}
	\left(\CGSRep(\CGSVect)-iJ\bigl(\CGSRep(\CGSVect)\bigr)\right)
\end{eqnarray*}
\end{definition}

With a little abuse of language we say that $\CGSCRep$ is
holomorphic if $\CGSCRep(\CGSVect)$ is a holomorphic vector field on $\CXMan$
for each $\CGSVect\in\CGSBase$.

With this notation the last proposition can be restated as follows.

\begin{theorem}\label{stm::CGSCommuting}
Let $(\CGSBase,\CGSRep,\CGSGradmap)$ be a complex gradient system.
Then for each $\CGSVect,\CGSVectB\in\CGSBase$ we have
\begin{eqnarray*}
\bigl[\CGSCRep(\CGSVect),\CGSCRep(\CGSVectB)]&=&0.
\end{eqnarray*}
\end{theorem}

}

\section{\label{section:Cauchy}A Cauchy problem}
{ 
Let $\CXMan$ be a complex manifold of complex dimension $\CXdim$.
Let $\CRMan\sset\CXMan$ be a $\CR-$submanifold 
of $\CXMan$ of type $(\CRdim,\CRcodim)$.

\begin{definition}
Let $\CGSBase$ be a real vector space.
A linear map $\CGSRep_0:\CGSBase\to\Gamma(\CRMan,\Tan\CRMan)$ is said to be
\emph{$\CR-$transverse} if for each $\CGSVect\in\CGSBase\setminus\{0\}$
and each $\pMan\in\CRMan$ we have
$J\bigl(\CGSRep_0(\CGSVect)(\pMan)\bigr)\not\in\Tan_\pMan\CRMan$.
\end{definition}

Given a vector field $\VField\in\Gamma(\CRMan,\Tan\CRMan)$ we denote
by $\Exp_\pMan(\VField)$ the exponential mapping associated to the
vector field $\VField$: if $\gamma(t)$ is the integral curve
of the vector field $\VField$ such that $\gamma(0)=\pMan$ then
$\Exp_\pMan(\VField)=\gamma(1)$.

Let $\CGSRep_0:\CGSBase\to\Gamma(\CRMan,\Tan\CRMan)$ be a linear map of real vector spaces.
The \emph{flow} associated to $\CGSRep_0$
is defined for $\pMan\in\CRMan$ and $\CGSVect\in\CGSBase$ by
\begin{eqnarray*}
	&&\CGSExp_{\CGSRep_0}(\pMan,\CGSVect)=
		\Exp_\pMan\bigl(\CGSRep_0(\CGSVect)\bigr)
		\end{eqnarray*}
$\CGSExp_{\CGSRep_0}$ is a smooth map which is well defined in
an open neighbourhood of
$\CRMan\times\{0\}$ in $\CRMan\times\CGSBase$.

Let $\CGSBaseC$ denote the complexification $\C\otimes_\R\CGSBase$ of the real vector space $\CGSBase$.
 We then say that the flow $\CGSExp_{\CGSRep_0}$ is
\emph{uniformly (real) analytic} if there exist
an open neighbourhood $\CGSCexpDom\sset\CRMan\times\CGSBaseC$
of $\CRMan\times\{0\}$ and a smooth function
\begin{eqnarray*}
	&&\CGSCExp_{\CGSRep_0}:\CGSCexpDom\to\CRMan
\end{eqnarray*}
which coincides with $\CGSExp_{\CGSRep_0}$
on $\CRMan\times\CGSBase$
and for each $\pMan\in\CRMan$ the map defined on the
open set
\begin{eqnarray*}
	&&\CGSCexpDom_\pMan=\bigl\{\CGSVect\in\CGSBaseC\mid
		(\pMan,\CGSVect)\in\CGSCexpDom\bigr\}
\end{eqnarray*}
by
\begin{eqnarray*}
	&&\CGSVect\mapsto\widetilde\CGSExp_{\CGSRep_0}(\pMan,\CGSVect)
\end{eqnarray*}
is holomorphic.

The map $\CGSCExp_{\CGSRep_0}$ will be called the \emph{complexification} of
the flow $\CGSExp_{\CGSRep_0}$.

Let $\CGSIExp_{\CGSRep_0}$ denote the restriction of
$\CGSCExp_{\CGSRep_0}$ to $\FIsoDom=\CGSCexpDom\cap\CRMan\times i\CGSBase$. As it is immediately seen, shrinking the domain $\CGSCexpDom$ if necessary, the map $\CGSIExp_{\CGSRep_0}$ is a diffeomorphism between
$\FIsoDom$ and $\CGSIExp_{\CGSRep_0}(\FIsoDom)$ if, and only if, the map $\CGSRep_0$ is $\CR-$transverse.

In this case we denote by
$\CGSGradmap_{\CGSRep_0}:\CGSIExp_{\CGSRep_0}(\FIsoDom)\to\CGSBase$
the unique map satisfying 
\begin{eqnarray*}
	\CGSGradmap_{\CGSRep_0}\bigl(\CGSIExp_{\CGSRep_0}(\pMan,i\CGSVect)\bigr)=\CGSVect.
\end{eqnarray*}
for each $\pMan\in\CRMan$ and each $\CGSVect\in\CGSBase$.

Observe that the map $\CGSGradmap_{\CGSRep_0}$
vanishes exactly on $\CRMan$, so we will refer to it as to the \emph{equation} of $\CRMan$ associated to $\CGSRep_0$.

We say that a complex gradient system
$(\CGSBase,\CGSRep,\CGSGradmap)$
\emph{extends} $\CGSRep_0$ if it is defined in a open
neighbourhood $\CGSCDom\sset\CXMan$ of $\CRMan$
and for every $\pMan\in\CRMan$ and $\CGSVect\in\CGSBase$
\begin{eqnarray*}
	\CGSRep(\CGSVect)(\pMan)=\CGSRep_0(\CGSVect)(\pMan).
\end{eqnarray*}
If $(\CGSBase,\CGSRep_1,\CGSGradmap)$, $(\CGSBase,\CGSRep_2,\CGSGradmap)$ are two extensions of $\CGSRep_0$ such that the for every $V\in \CGSBase$ the sections $\CGSRep_1(\CGSVect)$, $\CGSRep_2(\CGSVect)$ coincide on $N$ then we write ${\CGSRep_1}_{|\CGSCDom}={\CGSRep_2}_{|\CGSCDom}$.

\begin{theorem}\label{stm::CGSCauchy}
Let $\CXMan$ be a complex manifold of complex dimension $\CXdim$, $\CRMan\sset\CXMan$ a $\CR-$submanifold 
of $\CXMan$ of type $(\CRdim,\CRcodim)$. Let $\CGSBase$ be a real vector space and $\CGSRep_0:\CGSBase\to\Gamma(\CRMan,\Tan\CRMan)$ a $\CR-$transverse linear map such that the distribution $\CGSRDistr{\CGSRep_0}$ is integrable.
Assume that the associated flow $\CGSExp_{\CGSRep_0}$ is uniformly real analytic and let $\CGSGradmap_{\CGSRep_0}$ be the associated equation.

Then there exists an open neighbourhood $\CGSCDom\sset\CXMan$ of $\CRMan$
and an $\R$-linear map 
$$
\CGSRep:\CGSBase\to\Gamma(\CGSCDom,T\CXMan)
$$
such that $(\CGSBase,\CGSRep,\CGSGradmap_{\CGSRep_0})$
is a complex gradient system which extends $\CGSRep_0$.

The map $\CGSRep$ is unique in a neighbourhood of $\CRMan$, that is
if 
$$
\CGSRep_1:\CGSBase\to\Gamma(\CGSCDom_1,T\CXMan), \>\>\CGSRep_2:\CGSBase\to\Gamma(\CGSCDom_2,T\CXMan)
$$
are $\R$-linear maps such that $(\CGSBase,\CGSRep_1,\CGSGradmap_{\CGSRep_0})$
and $(\CGSBase,\CGSRep_2,\CGSGradmap_{\CGSRep_0})$ are
complex gradient systems which extend $\CGSRep_0$ then
${\CGSRep_1}_{|\CGSCDom}={\CGSRep_2}_{|\CGSCDom}$ for a suitable
open neighbourhood $\CGSCDom\sset\CGSCDom_1\cap\CGSCDom_2$ of $\CRMan$.
\end{theorem}

{
\def\FIso{{\CGSIExp_{\CGSRep_0}}}
\def\myJ{J}
\def\JOp#1{J\left({#1}\right)}
\def\hV#1{{\tilde\CGSVBase}_{#1}^0}
\def\vV#1{\frac{\partial}{\partial\ucomp_{#1}}}
\def\ISMan{S}
\def\matA{A}
\def\matP{P}
\def\matQ{Q}
\def\mata{a}
\def\matp{p}
\def\matq{q}
\def\ZNull{Z}
\proof  
It is not restrictive to assume that the map
$\CGSCExp_{\CGSRep_0}$ is a diffeomorphism between
$\FIsoDom=\CGSCexpDom\cap\CRMan\times i\CGSBase$ and $\FIso(\FIsoDom)$.

Fix a basis $\{\CGSVect_1,\ldots,\CGSVect_\CRcodim\}$ of $\CGSBase$
and set
\begin{eqnarray*}
	&&\CGSVBase_\ixa^0=\CGSRep(\CGSVect_\ixa),\quad\ixa=1,\ldots,\CRcodim.
\end{eqnarray*}

Let $\myJ$ denote the complex structure on $\Tan\FIsoDom$ induced by the pullback of 
of the complex structure on $\Tan\FIso(\FIsoDom)\sset\Tan\CXMan$
under the map $\FIso$.

It is not restrictive to identify the neighbourhood $\FIso(\FIsoDom)$ of
$\CRMan$ in $\CXMan$ with the domain $\FIsoDom\sset\CRMan\times i\CGSBase$.
We also identify $\CRMan\times i\CGSBase$ with $\CRMan\times\R^\CRcodim$
by
\begin{eqnarray*}
	&&\CRMan\times\R^\CRcodim\ni
	(\pMan,\ucomp_1,\ldots,\ucomp_\CRcodim)\mapsto
	(\pMan,i\ucomp_1\CGSVect_1+\cdots+i\ucomp_\CRcodim\CGSVect_\CRcodim)
	\in\CRMan\times i\CGSBase.
\end{eqnarray*}

Let
$$\CGSGradmap=(\ucomp_1,\ldots,\ucomp_\CRcodim):
	\CRMan\times\R^\CRcodim\to\R^\CRcodim\simeq\CGSBase
$$
be the projection on the second factor.

We will prove the existence of the required complex gradient system
showing that there exist vector fields
\begin{eqnarray*}
	&&\CGSVBase_\ixa\quad\ixa=1,\ldots,\CRcodim
\end{eqnarray*}
defined in a suitable neighbourhood $\CGSCDom$ of
$\CRMan\times\{0\}$ in $\FIsoDom$
such that
\begin{eqnarray*}
	&&\dd\ucomp_\ixa(\CGSVBase_\ixb)=0 \quad\ixa,\ixb=1,\ldots,\CRcodim
\end{eqnarray*}
and
\begin{eqnarray*}
	&&\dc\ucomp_\ixa(\CGSVBase_\ixb)=\delta_{\ixa\ixb} \quad\ixa,\ixb=1,\ldots,\CRcodim.
\end{eqnarray*}

Let $\hV\ixa$, $\ixa=1,\ldots,\CRcodim$, denote the vector fields on $\FIsoDom$
which coincide with $\CGSVBase_\ixa^0$ on $\CRMan\times\{0\}$ and are
invariant under the action $\CGSBase\times(\CRMan\times i \CGSBase)$ given by
\begin{eqnarray*}
	&&\bigl(\CGSVectB,(\pMan,\CGSVect)\bigr)\mapsto(\pMan,\CGSVectB+\CGSVect).
\end{eqnarray*}

Let also $\CGSDistr$ denote the distribution 
on $\Tan(\CRMan\times i\CGSBase)\approx\Tan(\CRMan\times\R^\CRcodim)$
generated by the vector fields
\begin{eqnarray*}
	&&\hV1,\ldots,\hV\CRcodim,\vV1,\ldots,\vV\CRcodim.
\end{eqnarray*}
$\CGSDistr$ is completely integrable and the maximal integral submanifolds of $\CGSDistr$ are of the form
$\ISMan\times\R^\CRcodim$ where $\ISMan$ is a maximal
integral submanifold of the distribution $\CGSRDistr\CGSRep$.

By construction, the intersection of each maximal integral submanifolds
of $\CGSDistr$ with the domain $\FIsoDom$ is
a complex submanifold of $\FIsoDom$ of complex dimension $\CRcodim$.
Moreover, for each $\pMan\in\CRMan$ and each $\ixa=1,\ldots,\CRcodim$ we have
\begin{eqnarray*}
	&&\myJ(\hV\ixa)(\pMan)=\vV\ixa(\pMan).
\end{eqnarray*}

Let $\matP=(\matp_{\ixa\ixb})$, $\matQ=(\matq_{\ixa\ixb})$ be the square matrices of
order $\CRcodim$ with entry smooth function on $\FIsoDom$ defined by
\begin{eqnarray*}
	&&\matp_{\ixa\ixb}=\myJ(\hV\ixb)(\ucomp_\ixa),\\
	&&\matq_{\ixa\ixb}=\JOp{\vV\ixb}(\ucomp_\ixa).
\end{eqnarray*}
Observe that for each $\pMan\in\CRMan$ the matrices
$\matP\bigl((\pMan,0)\bigr)$ and $\matQ\bigl((\pMan,0)\bigr)$
are respectively the identity matrix and the zero matrix of order $\CRcodim$.

Let $\CGSCDom$ be the open neighbourhood of $\CRMan\times\{0\}$
in $\FIsoDom$ defined by
\begin{eqnarray*}
	&&\CGSCDom=\Bigl\{(\pMan,\ucomp_1,\ldots,\ucomp_\CRcodim)\in\FIsoDom\mid
		\det\matP\bigl((\pMan,\ucomp_1,\ldots,\ucomp_\CRcodim)\bigr)
		\neq0\Bigr\}.
\end{eqnarray*}
Denote $\matA=(\mata_{\ixa\ixb})$ the matrix $\matP^{-1}\matQ$,
and set
\begin{eqnarray*}
	&&\CGSVBase_\ixa=-\JOp{\vV\ixa}+
		\sum_{\ixb=1}^\CRcodim\mata_{\ixb\ixa}\myJ(\hV\ixb)
		\quad\ixa=1,\ldots,\CRcodim.
\end{eqnarray*}
Then
\begin{eqnarray*}
	&&\myJ(\CGSVBase_\ixa)=\vV\ixa-
		\sum_{\ixb=1}^\CRcodim\mata_{\ixb\ixa}\hV\ixb
		\quad\ixa=1,\ldots,\CRcodim
\end{eqnarray*}
and, in view of the $\myJ$-invariance of the distribution $\CGSDistr$
it follows that 
\begin{eqnarray*}
	&&\CGSVBase_1,\ldots,\CGSVBase_\CRcodim,
	\myJ(\CGSVBase_1),\ldots,\myJ(\CGSVBase_\CRcodim)
\end{eqnarray*}
generate the distribution $\CGSDistr$ on $\CGSCDom$. It is easy to check that 
\begin{eqnarray*}
	&&\dd\ucomp_\ixa(\CGSVBase_\ixb)=0 \quad\ixa,\ixb=1,\ldots,\CRcodim
\end{eqnarray*}
and
\begin{eqnarray*}
	&&\dc\ucomp_\ixa(\CGSVBase_\ixb)=\delta_{\ixa\ixb} \quad\ixa,\ixb=1,\ldots,\CRcodim,
\end{eqnarray*}
as required.

In order to prove the uniqueness of the map $\CGSRep$ assume that the complex gradient systems $(\CGSBase,\CGSRep_1,\CGSGradmap_{\CGSRep_0})$, $(\CGSBase,\CGSRep_2,\CGSGradmap_{\CGSRep_0})$ extend $\CGSRep_0$
and set $\gamma=\CGSRep_1-\CGSRep_2$.

We are going to prove that, after shrinking $\CGSCDom$ if necessary, ${\gamma}_{|\CGSCDom}=0$ showing before that the complex distributions $\CGSCDistr{\CGSRep_1}$ and $\CGSCDistr{\CGSRep_2}$ associated respectively to $\CGSRep_1$ and $\CGSRep_2$ coincide near to $M\times\{0\}$.

By hypothesis the distribution $\CGSRDistr{\CGSRep_0}$ is integrable and
its maximal integral submanifold are real submanifolds of $\CRMan$
of (real) dimension $\CRcodim$. For every $\pMan\in\CRMan$ consider the maximal integral submanifolds $\ISMan_1$, $\ISMan_2$ through $\pMan$ of $\CGSCDistr{\CGSRep_1}$ and $\CGSCDistr{\CGSRep_2}$ respectively. Since $\CGSRep_1$ and $\CGSRep_2$ both extend $\CGSRep_0$ it
follows that 
$$
\ISMan^\mathbb{R}=\ISMan_1\cap\CRMan=\ISMan_2\cap\CRMan
$$
is the maximal integral (real) submanifold of (real) dimension $\CRcodim$
of the distribution $\CGSRDistr{\CGSRep_0}$ through $\pMan$.
In view of the hypothesis of $\CR-$transversality, $\ISMan^\mathbb{R}$ is a totally real
submanifold of $\ISMan_1$ and $\ISMan_2$.
It follows that $\ISMan_1=\ISMan_2$.

We have so proved that the maximal integral submanifolds of the
distributions  $\CGSCDistr{\CGSRep_1}$ and $\CGSCDistr{\CGSRep_2}$
which meet the submanifold $\CRMan$ are the same and consequently,
after shrinking $\CGSCDom$ if necessarily, it follows that the distributions
$\CGSCDistr{\CGSRep_1}$ and $\CGSCDistr{\CGSRep_2}$ coincide on $\CGSCDom$.

Let now $\CGSVect\in\CGSBase$ be an arbitrary vector.
Then, $\gamma(\CGSVect)\in\Gamma(\CGSCDom,\CGSCRTan_{\CGSRep_1})$
and the above argument shows that
$\gamma(\CGSVect)\in\Gamma(\CGSCDom,\CGSCDistr{\CGSRep_1})$.
Since
$\Tan\CXMan=\CGSRDistr{\CGSRep_1}\oplus\CXJ\CGSRDistr{\CGSRep_1}\oplus\CGSCRTan_{\CGSRep_1}$
it follows that ${\gamma}(V)_{|\CGSCDom}=0$ and this ends the proof $\CGSVect\in\CGSBase$ being arbitrary.

\qed
}

When $\CRcodim=\dim\CGSBase=1$ the result above is contained in
\cite[Theorem 3.1]{article:TV2ArxivContactMA}
where a stronger uniqueness result was obtained.
Namely, if $(\CGSVBase_1,\ucomp_1)$ and $(\CGSVBase_2,\ucomp_2)$
are two one dimensional calibrated foliation such that
$\CGSVBase_1$ and $\CGSVBase_2$ both extend $\CGSVBase_0$ along
the hypersurface $\CRMan$ then $\CGSVBase_1=\CGSVBase_2$
and $\ucomp_1=\ucomp_2$ in a neighbourhood of $\CRMan$.

Such a uniqueness result does not hold for a
general complex gradient system.
Indeed, consider $\CXMan=\C$, $\CRMan=\R$ and
$$\CGSVBase_0=\der[x].$$
Then our construction yields the gradient map
$$\CGSGradmap(z)=\CGSGradmap(x+iy)=-\Re(z)=-y$$
and the vector field
$$\CGSVBase=\der[x],$$
but also the pair $(\CGSVBase_1,\CGSGradmap_1)$ where
$$\CGSGradmap_1(z)=\CGSGradmap_1(x+iy)=e^{-y}-1$$
and
$$\CGSVBase_1=e^y\der[x],$$
is a complex gradient system which extends $\CGSVBase_0$.

In the case $\CRcodim=1$ the condition $[\CGSVBase,\CXJ\CGSVBase]=0$
which is in the definition of one dimensional calibrated foliation given in \cite{article:TV2ArxivContactMA}
ensures the uniqueness for the Cauchy problem
(see \cite[Theorem 3.1]{article:TV2ArxivContactMA}).
It is not clear which is the right condition (if any) to add
in order to guarantee the uniqueness also in this non commutative setting.

See also the examples given in the next section.

}

\section{\label{section:ExampleLieGroups}Lie Groups}
{ 
\def\GMetr{B}
\def\Adj{{\rm Ad}}

\def\MatTwoAlg[#1,#2]{
	\left(
		\begin{array}{cc}
			\displaystyle{#1}&\displaystyle{#2}\\
			0&0
		\end{array}
	\right)
}

\def\MatTwoGrp[#1,#2]{
	\left(
		\begin{array}{cc}
			\displaystyle{#1}&\displaystyle{#2}\\
			0&1
		\end{array}
	\right)
}

\def\MatThreeAlg[#1,#2,#3]{
	\left(
		\begin{array}{ccc}
			0&\displaystyle{#1}&\displaystyle{#3}\\
			0&0&\displaystyle{#2}\\
			0&0&0
		\end{array}
	\right)
}

\def\MatThreeGrp[#1,#2,#3]{
	\left(
		\begin{array}{ccc}
			1&\displaystyle{#1}&\displaystyle{#3}\\
			0&1&\displaystyle{#2}\\
			0&0&1
		\end{array}
	\right)
}

\def\rep[#1,#2]{\widetilde{#1}_{#2}}

\def\argA{{\theta_1}}
\def\factorA{\left(\frac{x_1}{y_1}-\frac{1}{\argA}\right)}

\def\elAlgA{\alpha}
\def\elAlgB{\beta}
\def\ElAlgA{A}
\def\ElAlgB{B}
\def\elGrpA{x}
\def\elGrpB{y}
\def\elCGrA{z}
\def\elCGrB{w}
\def\ElGrpA{U}
\def\ElGrpB{V}

Let $\CGSCGroup$ be a complex Lie group of (complex) dimension $\CRcodim$
which is the complexification of a real Lie group $\CGSGroup$;
$\CGSGroup$ is totally real submanifold of $\CGSCGroup$.
Let $\CGSCAlg$ and $\CGSAlg$ be the Lie algebras of
$\CGSCGroup$ and $\CGSGroup$ respectively.

We identify $\CGSAlg$ (resp. $\CGSCAlg$) with the tangent space
to $\CGSGroup$ ($\CGSCGroup$)
at the origin and for each
$\CGSAElem\in\CGSAlg$ (resp. $\CGSAElem\in\CGSCAlg$)
we denote
by $L_\CGSAElem$ the corresponding left invariant vector
field on $\CGSGroup$ (resp. $\CGSCGroup$).
The complexification of the flow associated to the map
$\CGSCAlg\ni\CGSVect\mapsto L_\CGSVect$ is the map
\begin{eqnarray*}
	\CGSGroup\times\CGSCAlg\ni
	(\CGSGElem,\CGSVect)\mapsto
	\CGSGElem\exp(\CGSVect)\in\CGSCGroup
\end{eqnarray*}
being $\exp$ the standard exponential map
$\exp:\CGSCAlg\to\CGSCGroup$.
Let denote by $(\CGSAlg,\CGSRep,\CGSGradmap)$
the complex gradient system which extends
$\CGSVect\mapsto L_\CGSVect$.
Then we have the identity
\begin{eqnarray*}
	\CGSGradmap\bigl(\CGSGElem\exp(-i\CGSVect)\bigr)=\CGSVect.
\end{eqnarray*}
	If $\CGSCGroup$ be a complex reductive Lie group and
$\CGSGroup$ is a compact real form  for $\CGSCGroup$
then we have the Cartan decomposition of $\CGSCGroup$ 
\begin{eqnarray*}
	&&\CGSGroup\times\CGSAlg\to\CGSCGroup\\
	&&(\CGSGElem,\CGSAElem)\mapsto\CGSGElem\exp(i\CGSAElem).
\end{eqnarray*}

In this case the algebra $\CGSAlg$ admits a definite metric $\GMetr$,
invariant under the adjoint representation $\Adj_\CGSGroup$ of $\CGSGroup$,
inducing an isomorphism between the Lie algebra $\CGSAlg$ and its
dual $\CGSAlgD$.
With this identification the gradient map $\CGSGradmap$ is (up to the sign)
the moment map associated to a symplectic action of $\CGSGroup$ on $\CGSCGroup$.
See e.g. \cite{article:HeinznerSchwarzCartanDecomp} for details.

This example explain our terminology ``\emph{complex gradient system}''.

We would like to point out that in general, as shown by the examples below,
the representation $\CGSRep$ of the complex gradient system
extending the left representation $\CGSVect\mapsto L_\CGSVect$
is not the restriction of the left representation of $\CGSCAlg$.

Let $\CGSCGroup$ be the  matrix Lie group of the matrices of the form
\begin{eqnarray*}
	&&\MatThreeGrp[z_1,z_2,z_3]
\end{eqnarray*}
with $z_\ixa=x_\ixa+iy_\ixa\in\mathbb{C}$, $\ixa=1,2,3$
and let $\CGSGroup$ be the corresponding group with $z_i\in\mathbb{R}$.

Then $\CGSAlg$ is given by the matrices of the form
\begin{eqnarray*}
	&&\MatThreeAlg[u_1,u_2,u_3]
\end{eqnarray*}
with $u_\ixa\in\mathbb{R}$, , $\ixa=1,\ldots,3$.

Put
\begin{eqnarray*}
	&&E_1=\MatThreeAlg[1,0,0],\ 
	E_2=\MatThreeAlg[0,1,0],\ 
	E_3=\MatThreeAlg[0,0,1].
\end{eqnarray*}
Then $E_1,E_2,E_3$ is a basis of $\CGSAlg$.
Denoting $L_\ixa=L_{E_\ixa}$ we then have
\begin{eqnarray*}
	L_1&=&\der[x_1],\\
	L_2&=&\der[x_2]+x_1\der[x_3]+y_1\der[y_3],\\
	L_3&=&\der[x_3],\\
	\CXJ L_1&=&\der[y_1],\\
	\CXJ L_2&=&\der[y_2]+x_1\der[y_3]-y_1\der[x_3],\\
	\CXJ L_3&=&\der[y_3]
\end{eqnarray*}
Some computation yields for the gradient map $\CGSGradmap$ the expression
\begin{eqnarray*}
	\CGSGradmap(z_1,z_2,z_3)=-y_1E_1-y_2E_2-(y_3+x_1y_2)E_3
\end{eqnarray*}
and the representation $\CGSRep$ is given by
\begin{eqnarray*}
	&&\CGSRep(E_1)=\rep[E,1]=\der[x_1]+y_2\der[y_3]=L_1+y_2J(L_3),\\
	&&\CGSRep(E_2)=\rep[E,2]=\der[x_2]+x_1\der[x_3]=L_2-y_1J(L_3),\\
	&&\CGSRep(E_3)=\rep[E,3]=\der[x_3]=L_3.
\end{eqnarray*}
Observe that
\begin{eqnarray*}
	&&\LieBracket{\rep[E,1]}{\rep[E,2]}=\rep[E,3],\  
	\LieBracket{\rep[E,1]}{\rep[E,3]}=
	\LieBracket{\rep[E,2]}{\rep[E,3]}=0,\\
	&&\LieBracket{\CXJ\rep[E,1]}{\CXJ\rep[E,2]}=\rep[E,3],\  
		\LieBracket{\CXJ\rep[E,1]}{\CXJ\rep[E,3]}=
		\LieBracket{\CXJ\rep[E,2]}{\CXJ\rep[E,3]}=0,\\
	&&\LieBracket{\rep[E,i]}{\CXJ\rep[E,j]} = 0\quad i,j=1,2,3.
\end{eqnarray*}
It follows that the representation
$\CGSRep:\CGSAlg\to\Gamma(\CGSCGroup,\Tan\CGSCGroup)$ is a
Lie algebra isomorphism and the gradient map $\CGSGradmap$ is a harmonic function.


Let now $\CGSCGroup$ be the  matrix Lie group of the matrices of the form
\begin{eqnarray*}
	&&\MatTwoGrp[z_1,z_2]
\end{eqnarray*}
with $z_1,z_2\in\mathbb{C}$, $z_1\neq0$,
and let $\CGSGroup$ be the corresponding group with $z_1,z_2\in\mathbb{R}$.

The Lie algebra $\CGSAlg$ of $\CGSGroup$ is given by the matrices of the form
\begin{eqnarray*}
	&&\MatTwoAlg[u_1,u_2]
\end{eqnarray*}
with $u_1,u_2\in\mathbb{R}$.
The matrices
\begin{eqnarray*}
	&&E_1=\MatTwoAlg[1,0],\ \
	E_2=\MatTwoAlg[0,1],\ 
\end{eqnarray*}
forms a basis of the Lie algebra $\CGSAlg$ which satisfies the relation
\begin{eqnarray*}
	&&\LieBracket{E_1}{E_2}=E_2.
\end{eqnarray*}
The corresponding left invariant vector fields on $\CGSCGroup$ are given by
\begin{eqnarray*}
	L_1&=&x_1\der[x_1]+y_1\der[y_1],\\
	L_2&=&x_1\der[x_2]+y_1\der[y_1],\\
	\CXJ L_1&=&-y_1\der[x_1]+x_1\der[y_1],\\
	\CXJ L_2&=&-y_1\der[x_2]+x_1\der[y_2].
\end{eqnarray*}
After some computations we obtain that the gradient map is given by
\begin{eqnarray*}
	\CGSGradmap(z_1,z_2)=-\argA E_1-\frac{y_2\argA}{y1}E_2
\end{eqnarray*}
where
\begin{eqnarray*}
	\argA=\arctan\frac{y_1}{x_1}
\end{eqnarray*}
and the representation $\CGSRep$ satisfies
\begin{eqnarray*}
	\CGSRep(E_1)=\rep[E,1]&=&x_1\der[x_1]+y_1\der[y_1]
		+y_2\factorA\der[x_2]
		+y_2\der[y_2],\\
	\CGSRep(E_2)=\rep[E,2]&=&\frac{y_1}{\argA}\der[x_2].
\end{eqnarray*}

Observe that
\begin{eqnarray*}
	&&\LieBracket{\rep[E,1]}{\rep[E,2]}=\LieBracket{\CXJ\rep[E,1]}{\CXJ\rep[E,2]}=\rep[E,2],
\end{eqnarray*}
and
\begin{eqnarray*}
	&&\LieBracket{\rep[E,1]}{\CXJ\rep[E,1]}=
		\frac{2y_2}{y_1}\factorA\rep[E,2],\\
	&&\LieBracket{\rep[E,1]}{\CXJ\rep[E,2]}=-\factorA\rep[E,2],\\
	&&\LieBracket{\rep[E,2]}{\CXJ\rep[E,2]}=0,
\end{eqnarray*}
namely the representation
$\CGSRep:\CGSAlg\to\Gamma(\CGSCGroup,\Tan\CGSCGroup)$ is a
Lie algebra isomorphism but the gradient map $\CGSGradmap$ is
not a harmonic function and
the vector fields $\rep[E,1],\rep[E,2],\CXJ\rep[E,1]$ and $\CXJ\rep[E,2]$
are not a basis of a Lie sub-algebra of $\Gamma(\CGSCGroup,\Tan\CGSCGroup)$.

}


\section{\label{section:HoloAbelian}The holomorphic abelian case}
{ 
\def\CGSDom{\Omega}
\def\Func{F}
\def\Diff{\varphi}
\def\Chart{U}

Let $(\CGSBase,\CGSRep,\CGSGradmap)$ be a complex gradient system.

With a little abuse of language we say that such a complex gradient system is
\emph{holomorphic} if $\CGSCRep(\CGSVect)$ is a holomorphic vector field on $\CXMan$
for each $\CGSVect\in\CGSBase$.

We also say that it is \emph{abelian} if
\begin{eqnarray*}
&&\bigl[\CGSCRep(\CGSVect),\bar{\CGSCRep(\CGSVectB)}\bigr]=0
\end{eqnarray*}
for each pair of vectors $\CGSVect,\CGSVectB\in\CGSBase$.
Such a condition is equivalent to say that
for each pair of vectors $\CGSVect,\CGSVectB\in\CGSBase$
one has
\begin{eqnarray*}
&&\bigl[\CGSRep(\CGSVect),\CGSRep(\CGSVectB)\bigr]=
	\bigl[\CGSRep(\CGSVect),\CXJ\CGSRep(\CGSVectB)\bigr]=
	\bigl[\CXJ\CGSRep(\CGSVect)),\CXJ\CGSRep(\CGSVectB)\bigr]=0
\end{eqnarray*}
Consider now a domain $\CGSDom\sset\CC^\CRdim$ and
let $\Func:\CGSDom\to\RR^\CRcodim$ be a smooth function. We associate to $\Func$ a complex gradient system as follows.

Set $\CXMan=\CGSDom\times\CC^\CRcodim$, $\CGSBase=\RR^\CRcodim$
and define
\begin{eqnarray*}
\CGSGradmap(z,w)=\Func(x,y)-u,
\end{eqnarray*}
where $z=x+iy$ and $w=t+iu$ with $x,y\in\RR^\CRdim$ and $t,u\in\RR^\CRcodim$.
Finally consider the linear map
$\CGSRep:\RR^\CRcodim\to\Gamma(\CXMan,\Tan\CXMan)$ characterized 
by the conditions
\begin{eqnarray*}
&&\CGSRep(e_\ixa)=\der[t_\ixa]\quad\ixa=1,\ldots,\CRcodim,
\end{eqnarray*}
where $e_1,\ldots,e_\CRcodim$ is the canonical basis of $\RR^\CRcodim$.

It is easy to show that this complex gradient system is holomorphic and abelian and the aim of the next theorem is to prove that it is the local model of any holomorphic abelian complex gradient system. Namely the following is true
\begin{theorem}\label{stm::HoloAbelian}
Let $\CXMan$ be a complex manifold of complex dimension $\CXdim$.
Let $(\RR^\CRcodim,\CGSRep,\CGSGradmap)$ be a
holomorphic abelian complex gradient system on $\CXMan$.
Then for each point $\pMan$ there exist a complex coordinate system
\begin{eqnarray*}
z=(z_1,\ldots,z_n), \>\>w=(w_1,\ldots,w_k)
\end{eqnarray*}
$z_\iya=x_\iya+iy_\iya, \iya=1,\ldots,\CRdim$, $w_\ixa=t_\ixa+iu_\ixa, \ixa=1,\ldots,\CRcodim$, $x=(x_1,\ldots,x_\CRdim)$, $y=(y_1,\ldots,y_\CRdim)$,  $u=(u_1,\ldots,u_\CRcodim)$ and a smooth (vector) function $\Func$
depending only on $x$ and $y$ such that
\begin{eqnarray*}
&&\CGSRep(e_\ixa)=\der[t_\ixa]\quad\ixa=1,\ldots,\CRcodim,
\end{eqnarray*}
\begin{eqnarray*}
&&\CGSGradmap(z,w)=\Func(x,y)-u.
\end{eqnarray*}

\end{theorem}
\proof
Let 
$$
g_t^1,\ldots,g_t^\CRcodim, h_t^1,\ldots,h_t^\CRcodim
$$
be the (local) one parameter group of transformation of $\CXMan$ generated by the vector fields
$$
\CGSVBase_1,\ldots,\CGSVBase_\CRcodim, J(\CGSVBase_1),\ldots,J(\CGSVBase_\CRcodim).
$$
By hypotheses the Lie brackets between all pair of vector fields among
$\CGSVBase_1,\ldots,\CGSVBase_\CRcodim$ and
$\CXJ\CGSVBase_1,\ldots,\CXJ\CGSVBase_\CRcodim$ are zero
and hence the transformations
$g_t^1,\ldots,g_t^\CRcodim$ and $h_t^1,\ldots,h_t^\CRcodim$
commute each other.

Let $\pMan\in\CXMan$ be fixed and let $z_1,\ldots,z_\CXdim$ be a complex coordinates system around $\pMan$,
where $z_\iya=x_\iya+iy_\iya, \iya=1,\ldots,\CXdim$.

After reordering the coordinates we may suppose that
\begin{eqnarray*}
&&\der[x_1],\der[y_1],\ldots,
	\der[x_\CRdim],\der[y_\CRdim],
	\CGSVBase_1,\CXJ\CGSVBase_1\ldots,
	\CGSVBase_\CRcodim,\CXJ\CGSVBase_\CRcodim
\end{eqnarray*}
generates the tangent space to $\CXMan$ at each point
in a suitable neighbourhood of $\pMan$.

For $\ixa=1,\ldots,\CRcodim$ set $w_\ixa=t_\ixa+iu_\ixa$ and define
\begin{eqnarray*}
&&G^\ixa_{w_\ixa}=g^\ixa_{t_\ixa}\circ h^\ixa_{u_\ixa},
\end{eqnarray*}
Then the map
\begin{eqnarray*}
&&(z_1,\ldots,z_\CRdim,w_1,\ldots,w_\CRcodim)\mapsto
G^1_{w_1}\circ\cdots\circ G^\CRcodim_{w_\CRcodim}
(z_1,\ldots,z_\CRdim,0,\ldots,0),
\end{eqnarray*}
is a diffeomorphism $\Diff$ between an open set of
$\CC^\CXdim$ and a suitable neighbourhood $\Chart$
of $\pMan$ in $\CXMan$, that is
\begin{eqnarray*}
&&x_1,y_1,\ldots,x_\CRdim,y_\CRdim,t_1,u_1,\ldots,t_\CRcodim,u_\CRcodim
\end{eqnarray*}
is a real coordinate system on $\Chart$.

Since the maps $G^1_{w_1}\ldots G^\CRcodim_{w_\CRcodim}$
commute each other it follows that with respect to such a
coordinate system we have
\begin{eqnarray*}
&&\CGSVBase_\ixa=\CGSRep(e_\ixa)=\der[t_\ixa]
\end{eqnarray*}
for $\ixa=1,\ldots,\CRcodim.$ 

We now prove 
\begin{eqnarray*}
&&z_1,\ldots,z_\CRdim,w_1,\ldots,w_\CRcodim
\end{eqnarray*}
are complex coordinates on $\Chart$, showing that the diffeomorphism $\Diff$ is in fact a biholomorphism. 

Since, by hypotheses,  $(\RR^\CRcodim,\CGSRep,\CGSGradmap)$ is a holomorphic abelian complex gradient system,
it follows that $G^1_{w_1}\ldots G^\CRcodim_{w_\CRcodim}$
are holomorphic local diffeomorphisms and for fixed $w_1,\ldots,w_\CRcodim$ the map
\begin{eqnarray*}
&&\Diff(z_1,\ldots,z_\CRdim,w_1,\ldots,w_\CRcodim)
\end{eqnarray*}
is holomorphic with respect to the variables $z_1,\ldots,z_\CRdim$. Moreover, for $\ixa=1,\dots,\CRcodim$ the maps $g_{t_\ixa}$ and $h_{u_\ixa}$
commute and hence the map $w_\ixa\mapsto G^\ixa_{w_\ixa}(\cdot)$
is holomorphic with respect to $w_\ixa$.
On the other hand, since the maps $G^1_{w_1}\ldots G^\CRcodim_{w_\CRcodim}$
commute each other, the map $\Diff$ is holomorphic with respect to the variable $w_\ixa$, $\ixa=1,\dots,\CRcodim$,
when the variables $z_1,\ldots,z_\CRdim$ and $w_1,\ldots,w_{\ixa-1},$
$w_{\ixa+1}\ldots,w_\CRcodim$
are fixed.

Thus the map $\Diff$ is separately holomorphic in each variable
and hence is holomorphic.

Finally let 
$\CGSGradmap=(\CGSGradmap_1,\ldots,\CGSGradmap_\CRcodim):\CXMan\to\RR^\CRcodim$
be the gradient map
and consider the map $F=(F_1,\ldots,F_\CRcodim):\Chart\to\RR^\CRcodim$
defined by
\begin{eqnarray*}
&&F_\ixa(z,w)=F_\ixa(x,y,t,u)=\CGSGradmap_\ixa(z,w)+u_\ixa
	\quad\ixa=1,\ldots,\CRcodim.
\end{eqnarray*}
We end the proof showing that the map $F$ does not depend on the variables $t$ and $u$.
Indeed, for $\ixa,\ixb=1,\dots,\CRcodim$, we have
\begin{eqnarray*}
&&\pder{F_\ixa}{t_{\ixb}}=\CGSVBase_\ixb(\CGSGradmap_\ixa)=0
\end{eqnarray*}
and
\begin{eqnarray*}
&&\pder{F_\ixa}{u_{\ixb}}=J(\CGSVBase_\ixb)(\CGSGradmap_\ixa)+\delta_{\ixa\ixb}=
-\delta_{\ixa\ixb}+\delta_{\ixa\ixb}=0.
\end{eqnarray*}{}
\qed

}

\nocite{book:DusaSalamonIntro}
\nocite{article:HeinznerSchwarzCartanDecomp}
\nocite{book:AnaCannasDaSilva}


\begin{thebibliography}{1}

\bibitem{book:AnaCannasDaSilva}
A.~C. da~Silva.
\newblock {\em {Lectures on Symplectic Geometry}}, volume 1764 of {\em Lecture
  Notes in Mathematics}.
\newblock Springer-Verlag, 2001.

\bibitem{book:DragominTomassini}
S.~Dragomir and G.~Tomassini.
\newblock {\em {Differentiable Geometry and Analysis on CR Manifolds}}.
\newblock Birkhauser, 2006.

\bibitem{article:DuchampKalka}
T.~Duchamp and M.~Kalka.
\newblock {Singular Monge-Amp\`ere foliations}.
\newblock {\em Math. Ann.}, 325:187--209, 2003.

\bibitem{article:HeinznerSchuetzdellerConvexGradMap}
P.~Heinzner and P.~Schuetzdeller.
\newblock {Convexity properties of gradient maps.}
\newblock {\em {arXiv:0710.1152v1}}, pages 1--16, 2007.

\bibitem{article:HeinznerSchwarzCartanDecomp}
P.~Heinzner and G.~Schwarz.
\newblock {The Cartan decomposition of the moment map.}
\newblock {\em Math. Ann.}, 337:197--232, 2007.

\bibitem{book:DusaSalamonIntro}
D.~McDuff and D.~Salomon.
\newblock {\em {Introduction to Symplectic Topology}}.
\newblock Oxford Mathematical Monographs. Oxford University Press, 1995.

\bibitem{article:TV2ArxivContactMA}
G.~Tomassini and S.~Venturini.
\newblock {Contact geometry of one dimensional holomorphic foliations}.
\newblock {\em arXiv:0907.5082v1}, pages 1--15, 2009.
\newblock To appear in Indiana J. of Math.

\bibitem{article:TV3ArxivAdaptedTube}
G.~Tomassini and S.~Venturini.
\newblock {Adapted complex tubes on the symplectization of pseudo-Hermitian
  manifolds}.
\newblock {\em {arXiv:1002.4558}}, pages 1--6, 2010.

\bibitem{article:TV3AdaptedTube}
G.~Tomassini and S.~Venturini.
\newblock {Adapted complex tubes on the symplectization of pseudo-Hermitian
  manifolds}.
\newblock {\em Arkive Math.}, 96:77--83, 2011.

\end{thebibliography}
\end{document}